\documentclass[12pt]{amsart}
\usepackage{amsfonts}
\usepackage{amsmath,amssymb,amsthm}
\usepackage{amstext}
\usepackage{verbatim}
\input{epsf}
\newtheorem{thm}{Theorem}[section]

\newtheorem{prop}[thm]{Proposition}

\theoremstyle{definition}

\theoremstyle{remark}
\newtheorem{rem}[thm]{Remark}

\begin{document}

\title[Direct and indirect methods of proof. The Lehmus-Steiner theorem.]
{Direct and indirect methods of proof. The Lehmus-Steiner's theorem}
\author{Vesselka Mihova and  Julia Ninova}%
\address{Faculty of Mathematics and Informatics, University of Sofia
e-mail: mihova@fmi.uni-sofia.bg}%
\address{Faculty of Mathematics and Informatics, University of Sofia
e-mail: julianinova@hotmail.com}


\maketitle \thispagestyle{empty}

\begin{abstract}
We describe and discuss  different methods of proof of a
given statement and illustrate by logical models
the essence of specific types of proofs, especially of direct and indirect methods of proof.

Direct proofs of Lehmus-Steiner's Theorem are proposed.
\end{abstract}

\section{Introduction}
In most fields of study, knowledge is acquired by way of
observations, by reasoning about the results of observations and
by studying the observations, methods and theories of other fields
and practices.

Ancient Egyptian, Babylonian and Chinese mathematics consisted of
rules for measuring land, computing taxes, predicting eclipses,
solving equations and so on.

The ancient Greeks found that in arithmetic and geometry it was
possible to \emph{prove} that observation results are true. They
found that some truths in mathematics were \emph{obvious} and that
many of the others could be shown to follow logically from the
obvious ones.

On the other hand, Physics, Biology, Economics and other sciences
\emph{discover general truths} relying on observations. Besides,
not any general truth can be proved to be true - it can only be
tested for contradictions and inconsistencies. If a scientific
theory is accepted because observations have agreed with it, there
is in principle small doubt that a new observation will not agree
with the theory, even if all previous observations have agreed
with that theory. However, if a result is proved thoroughly and
correctly, that cannot happen.

Under what conditions can we be sure that the steps in our
investigations are correct? Are we really sure that what seems to
be obvious to us is in fact true? Can we expect \emph{all}
mathematical truths to follow from the obvious ones? These
questions are not easily answered.

Disputes and mistakes about what is obvious could be avoided by
laying down certain basic notions, relations and statements, called
\emph{axioms} (postulates assumed true, but unprovable) for each
branch of mathematics, and agreeing that proofs of assertions must
be derived from these. To axiomatize a system of knowledge means to
show that its claims can be derived from a small, well-understood
set of axioms (see also \cite{D}).

The
axiomatic system is subordinated to some conditions.
\begin{itemize}
\item[-] The system must be consistent, to lack contradiction, i. e. the ability to derive both a statement and
its negation from the system's axioms.

Consistency is a necessary requirement for the system. \vskip 1mm

\item[-] Each axiom has to be independent, i. e. not a theorem that can be derived from other axioms in the system.

However, independence is not a necessary requirement for the
system. \vskip 1mm

\item[-] The system can be complete, i. e. for every statement, either itself or its negation is derivable.
\end{itemize}

There is no longer an assumption that axioms are \emph{true in any
sense}; this allows parallel mathematical theories to be built on
alternate sets of axioms (for instance  \emph{Axiomatic set
theory, Number theory}).  \emph{Euclidean} and \emph{Non-Euclidean geometry}
have a common basic set of axioms;  the differences between these important geometries
are based on their alternate axioms of parallel lines.

Another way to avoid mistakes about what is obvious in mathematics
could be the use of rules of inference with purely formal content.

In mathematical logic a propositional calculus (also called sentential calculus or sentential logic)
is a formal system in which formulas of a formal language may be interpreted to represent propositions.

In  \cite {MN2, MN, MN1, MN3} we explain methods, based on logical
laws, for composition and proof of equivalent and inverse
 problems.

In \cite {MN} we discuss a way of generating groups of equivalent
problems. The method we propound is based on the logical
equivalences
$$p\wedge \neg q \rightarrow r\;\Leftrightarrow\;p\wedge \neg r\rightarrow q
\;\Leftrightarrow\; p\rightarrow q\vee r,$$
where $\,p, q, r\,$ are statements.

Using the sentential logic in \cite {MN1} and  \cite {MN3} we
propose a new problems composing technology as an interpretation of
specific logical models.
Clarifying and using the logical equivalence (see also  \cite {S})
$$(t\wedge p\rightarrow r)\,\wedge\,(t\wedge q\rightarrow r)
\;\Leftrightarrow\; t\wedge (p\vee q)\rightarrow r,\leqno(*)$$ we
give  an algorithm for composition of \emph{inverse} problems with a given
logical structure that is based on the steps below.
\begin{itemize}
\item[-] Formulating  and proving  \emph{generating} problems with logical structures
of the statements as those at the left hand side of (*).

\item[-] Formulating a problem with logical structure $\;t\wedge (p\vee q)\rightarrow r\,$
of the statement.

\item[-] Formulating  and proving the \emph{inverse} problem
with logical structure $\,t\wedge r \rightarrow  p\vee q$.
\end{itemize}

In \cite{MN2}, besides the generalization of criteria \emph{A} and \emph{D} for congruence of triangles,
we also illustrate the above algorithm  by suitable groups of examples.
\vskip 2mm

In section 2 of the present paper we describe and discuss different methods of
proof of implicative statements and illustrate by logical models
the essence of specific types of proofs, especially of direct and indirect
proofs.
\vskip 2mm

In section 3 of the paper we  propose \emph{direct proofs} of
Lehmus-Steiner's Theorem that differ from any we have come across.
\vskip 1mm

Our investigations in this field are appropriate for training of mathematics students and teachers.

\vskip 6mm
\section{Types of Proofs.}

Both discovery and proof are integral parts of problem solving.
The \emph{discovery} is thinking of possible solutions, and the
\emph{proving} ensures that the proposed solution actually solves
the problem.

\emph{Proofs} are logical descriptions of deductive reasoning and are
distinguished from inductive or empirical arguments; a proof must
demonstrate that a statement is always true (occasionally by
listing all possible cases and showing that it holds in each).

An unproven statement that is believed true is known as a
\emph{conjecture}.

The \emph{objects} of proofs are premises, conclusions, axioms,
theorems (propositions derived earlier from axioms), definitions and evidence from the real world.

The abilities (techniques) to have a working knowledge of these objects include
\begin{itemize}
\item[-] Rules of inference: simple valid argument forms.
They may be divided into \emph{basic rules}, which are fundamental
to logic and cannot be eliminated without losing the ability to
express some valid argument forms, and \emph{derived rules}, which
can be proven by the basic rules.

 To sum up, the rules of inference are logical rules which allow the deduction of conclusions from premises.
\vskip 2mm

\item[-] Laws of logical equivalence.
\end{itemize}

Different methods of proof combine these objects and techniques in different ways to
create valid arguments.
\vskip 2mm

According to \emph{Euclid} a precise  proof of a given statement
has the following structure:

- \emph{Premises}: These include given axioms and theorems, true
statements, strict restrictions for the validity of the given
statement, chosen suitable denotations. (\emph{It is given}...)

- \emph{Statement}: Strict formulation of the submitted statement.
(\emph{It is to be proved that}...)

- \emph{Proof}: Establishing the truth of the submitted statement
using premises, conclusions, rules of inference and logical laws.
\vskip 3mm

Let now $P$ and $Q$ be statements. In order to establish the truth of
the implication $P\rightarrow Q$, we discuss different methods of
proof. Occasionally, it may be helpful first to rephrase certain
 statements, to clarify that they are really formulated in an implicative form.

If ``not''  is put in front of a statement $P$, it negates the
statement. $\neg P$ is sometimes called the \emph{negation} (or
\emph{contradictory}) of $P$. For any statement $P$ either $P$ or
$\neg P$ is true and the other is false.
\vskip 4mm

\textbf{Formal Proofs.} The concept of a proof is formalized in
the field of mathematical logic. Purely formal proofs, written in
symbolic language instead of natural language, are considered in
proof theory. A formal proof is defined as a sequence of formulas in
a formal language, in which each formula is a logical consequence
of preceding formulas.

In a formal proof the statements $P$ and $Q$ aren't necessarily
related comprehensively to each other. Only the structure of the statements and the logical rules
that allow the deduction of conclusions from premises are
important.

Hence, to prove formally that an argument $Q$ is valid or the conclusion
follows logically from the hypotheses $P$, we have to

-  assume the hypotheses $P$ are true,

- use the formal rules of inference and logical equivalences to
determine that the conclusion $Q$ is true.

The following logical equivalences illustrate a formal proof:
$$\neg{(P\rightarrow Q)}\Leftrightarrow\neg{(\neg{P}\vee Q)}\Leftrightarrow
\neg{(\neg{P})}\wedge \neg Q\Leftrightarrow P\wedge \neg Q.$$
\vskip 4mm

\textbf{Vacuous proof.} A vacuous proof of an implication happens when the
hypothesis of the implication is always false, i. e.
if we know one of the hypotheses in $P$ is
false then $P\rightarrow Q$ is vacuously true.

For instance, in the implication $(P\wedge\neg P)\rightarrow Q$ the hypotheses
form a contradiction. Hence, $Q$ follows from the hypotheses
vacuously.
\vskip 4mm

\textbf{Trivial proofs.} An implication is trivially true when its conclusion is always true.
 Consider an implication $P\rightarrow Q$.
 If it can be shown (independently of $P$) that $Q$ is true, then the
implication is always true.

The form of the trivial proof
$\;Q \rightarrow (P \rightarrow Q)$ is, in fact, a tautology.
\vskip 4mm

\textbf{Proofs of equivalences.} For equivalence proofs or proofs of
statements of the form $P$ \emph{if and only if} $Q$  there are two
methods.
\begin{itemize}
\item[-] Truth table.

\item[-] Using direct or indirect methods and the equivalence
$$(P\leftrightarrow Q)\;\Leftrightarrow\;(P\rightarrow
Q\,\wedge\, Q\rightarrow P).$$

Thus, the proposition $P$ \emph{if and only if} $Q$ can be
proved if both the implication $P\rightarrow Q$ and the
implication $Q\rightarrow P$ are proved. This is the definition of the \emph{biconditional} statement.
\end{itemize}
\vskip 4mm

\textbf{Proof by cases.} If the hypothesis $P$ can be separated
into cases $p_1\vee p_2\vee ... \vee p_k$,  each of the
propositions $p_1\rightarrow Q$, $p_2 \rightarrow Q$, . . . , $p_k
\rightarrow Q$, is to be proved separately. A statement $P\rightarrow Q$ is true if
all possible cases are true.

The logical equivalences in this case are  (see also \cite {S}, p. 81)
$$p_1\rightarrow Q \wedge p_2\rightarrow Q \wedge...\wedge p_k\rightarrow Q
\;\Leftrightarrow\; p_1\vee p_2\vee...\vee p_k \rightarrow Q\;\Leftrightarrow\; P\rightarrow Q.$$

Different methods may be used to prove the different cases.
\vskip 4mm

\textbf{Direct proof.} In mathematics and logic, a direct proof is a
way of showing the truth or falsehood of a given statement by a
straightforward combination of established facts, usually existing
lemmas and theorems.

The methods of proof of these established facts,
lemmas, propositions and theorems are of \emph{no importance}.
Their truth or falsehood are to be accepted \emph{without any effort}.

However, it is exceptionally important that the actual proof of the given statement consists of
straightforward combinations of these facts  \emph{without making any further assumptions}.
\vskip 2mm

Thus, to prove an implication $P\rightarrow Q$ directly, we assume that
statement $P$ holds and try to deduce that statement $Q$
\emph{must} follow.
\vskip 2mm

The structure of the direct proof is:
\begin{itemize}
\item[-] \emph{Given} - a statement of the form $P\rightarrow Q$.

\item[-] \emph{Assumption} - the hypotheses in $P$ are true.

\item[-] \emph{Proof} - using the rules of inference, axioms, theorems and any
logical equivalences to establish in a straightforward way the truth of the conclusion $Q$.
\end{itemize}
\vskip 4mm

\textbf{Indirect proof.} It is often very difficult to give a
direct proof to $P \rightarrow Q$. The connection between $P$ and
$Q$ might not be suitable to this approach.

Indirect proof is a type of proof in which a statement to be
proved is assumed false and if the assumption leads to an
impossibility, then the statement assumed false has been proved to
be true.

There are four possible implications we can derive from the
implication $P\rightarrow Q$, namely
\begin{itemize}
\item[-] \emph{Conversion} (the \emph{converse}): $Q\rightarrow P$,

\item[-] \emph{Inversion} (the \emph{inverse}): $\neg P\rightarrow \neg Q$,

\item[-] \emph{Negation}: $\neg(P\rightarrow Q)$,

\item[-] \emph{Contraposition} (the \emph{opposite}, \emph{contrapositive}): $\neg Q\rightarrow \neg P$.
\end{itemize}

The implications $P\rightarrow Q$ and $\neg Q\rightarrow\neg P$
are logically equivalent.

The implications $Q\rightarrow P$ and $\neg P\rightarrow \neg Q$
are logically equivalent too, but they are not equivalent to the
implication $P\rightarrow Q$.

\vskip 3mm

The two most common indirect methods of proof are called
\emph{Proof by Contraposition} and \emph{Proof by Contradiction}.
These methods of indirect proof differ from each other in the assumptions we do
as premisses.
\vskip 4mm

\begin{itemize}
\item[] \textbf{Proof by Contraposition.}
In logic, contraposition is a law that says that a conditional
statement is logically equivalent to its contrapositive. This is
often called the \emph{law of contrapositive}, or the \emph{modus
tollens} (\emph{ denying the consequent}) rule of inference.
\vskip 2mm

The structure of this indirect proof is:
\begin{itemize}
\item[-] We consider an implication $P\rightarrow Q$.

\item[-] Its contrapositive (opposite) $\neg Q\rightarrow\neg P$ is logically equivalent to
the original implication, i. e.
$$\neg Q\rightarrow\neg P \;\Leftrightarrow \;P\rightarrow Q.$$

\item[-] We prove that\emph{ if $\neg Q$ is true} (the assumption), then $\neg P$ is true.
\end{itemize}
\vskip 2mm

Therefore, a proof by contraposition is a \emph{direct} proof of the
contrapositive.
\vskip 3mm

The proof of Lehmus-Steiner's Theorem in \cite{St} is
an illustration of a proof by contraposition.
\vskip 4mm

\item[] \textbf{Proof by contradiction.}
In logic, proof by contradiction is a form of proof, and more
specifically a form of indirect proof, that establishes the truth
or validity of a proposition by showing that the proposition's
being false would imply a contradiction. Proof by contradiction is
also known as \emph{indirect} proof, \emph{apagogical argument}, proof by
assuming the opposite, and \emph{reductio ad impossibility}. It is
a particular kind of the more general form of argument known as
\emph{reductio ad absurdum}.

We assume the proposition $P\rightarrow Q$ is false by assuming the negation of the
conclusion $Q$ and the premise $P$ are true, and then using $P
\wedge \neg Q$ to derive a contradiction.
\vskip 2mm

Hence, the structure of
this indirect proof is:
\begin{itemize}
\item[-] We use the equivalence $(P \rightarrow Q) \Leftrightarrow (\neg P \vee Q)$.

\item[-] The negation of the last disjunction is $P \wedge \neg Q$, i. e.\\
 $\neg (P \rightarrow Q)\Leftrightarrow (P \wedge \neg Q)$.

\item[-] To prove the original implication $P \rightarrow Q$,
we show that \emph{if its negation $P \wedge \neg Q$ is true} (the assumption), then this leads to a
contradiction.
\end{itemize}
\vskip 2mm

In other words, to prove the implication $P \rightarrow Q$ by
contradiction, we assume the hypothesis $P$ and the negation of
the conclusion $\neg Q$ both hold and show that this is a
contradiction (see also \cite {S}, p. 188).
\vskip 2mm

A logical base of this method are equivalences of the form
$$\begin{array}{lll}
P\rightarrow Q &\Leftrightarrow \neg Q\wedge P \rightarrow \neg P &
\Leftrightarrow \neg (P\rightarrow Q)\rightarrow\neg P;\\
[2mm]
P\rightarrow Q &\Leftrightarrow \neg Q\wedge P \rightarrow Q &
\Leftrightarrow \neg (P\rightarrow Q)\rightarrow Q.
\end{array}$$

Let now $T$ be a valid theorem, statement, axiom or definition of a notion in
the corresponding system of knowledge. The following
equivalences can also be logical base of a Proof by Contradiction of
the implication $\,P\rightarrow Q$.
$$P\rightarrow Q \;\Leftrightarrow \neg Q\wedge P \rightarrow \neg T \;
\Leftrightarrow \neg (P\rightarrow Q)\rightarrow\neg T.$$

The theoretical base of this method of proof is the \emph{law of excluded middle}
(or the \emph{principle of excluded middle}). It states that for any proposition,
either that proposition is true, or its negation is true.
The law is also known as the \emph{law} (or \emph{principle}) \emph{of the excluded third}.
\end{itemize}
\vskip 2mm

Examples of indirect proofs of Lehmus-Steiner's Theorem are given in \cite {K}.
\vskip 2mm

There exist also examples of indirect proofs of implications $P \rightarrow Q$ in which the statement $\neg Q$
can be separated into cases $q_1\vee q_2\vee ... \vee q_k$, $k\geq 2,\, k\in \mathbb{N}$.
In such a case each of the propositions $P\rightarrow q_1$, $P \rightarrow q_2$, . . . , $P
\rightarrow q_k$ is to be proved separately to be false. If moreover the premise $P$ is true it follows that
all the statements $q_i,\, i=1,...,k,$ are false and the conclusion $Q$ is true,
i. e.
$$\neg(\neg Q)\,\Leftrightarrow \,\neg(q_1\vee q_2\vee ... \vee q_k)\;\Leftrightarrow\;
\neg q_1\wedge\neg q_2\wedge...\wedge\neg q_k\,\Leftrightarrow\, Q.$$

The logical equivalences in this case are  (see also \cite {S}, p. 81)
$$P\rightarrow \neg q_1\,\wedge\, P\rightarrow \neg q_2\,\wedge
...\wedge \,P\rightarrow \neg q_k\;\Leftrightarrow\;
P\,\rightarrow\,\neg q_1\wedge\neg q_2\wedge...\wedge\neg q_k
\; \Leftrightarrow\; P\,\rightarrow\, Q.$$

The indirect proof of Lehmus-Steiner's theorem given in \cite {HK}
has in fact logical structure as the described above although this is not
mentioned by the authors.
\vskip 4mm

\textbf{Proof by construction.} In mathematics, a
\emph{constructive} proof is a method of proof that demonstrates
the existence of a mathematical object by creating or providing a
method for creating the object.

In other words, \emph{proof by construction} (proof by example) is
the construction of a concrete example with a property to show
that something having that property exists.

A simple constructive proof of Lehmus-Steiner's Theorem is given in \cite{W}.
\vskip 4mm

\textbf{Nonconstructive proof.} A nonconstructive proof
establishes that a mathematical object with a certain property
exists without explaining how such an object can be found.
This often takes the form of a proof by contradiction in which the
nonexistence of the object is proven to be impossible.
\vskip 4mm

\textbf{Proof by counterexamples.} We can \emph{disprove}
something by showing a single counter example, i. e. one finds an
example to show that something is not true.

However, we cannot prove something by example.
\vskip 4mm

\textbf{Mathematical induction.}  In proof by mathematical
induction, a single \emph{base case} is proved, and an
\emph{induction rule} is proved, which establishes that a certain
case implies the next case. Applying the induction rule
repeatedly, starting from the independently proved base case,
proves many, often infinitely many, other cases. Since the base
case is true, the infinity of other cases must also be true, even
if all of them cannot be proved directly because of their infinite
number.

The mathematical induction is a method of mathematical proof
\emph{typically} used to establish a given statement for all
natural numbers. It is a form of direct proof and it is done in
three steps.

Let $\mathbb{N} = \{1,2,3,4,...\}$ be the set of natural numbers,
and $P(n)$ be a mathematical statement involving the natural
number $n\geq k$, $k, n\in \mathbb{N}$, $k$ suitably fixed.

- The first step, known as the \emph{base step}, is to prove the
given statement for the first possible (admissible) natural number
$k$, i.e. to show that $P(k)$ is true for $n=k$. \vskip 1mm

- The second step, known as the \emph{inductive hypothesis}, is to
assume that for a natural number $i \geq k$ the statement
$P(i),\; i\in \mathbb{N}$ is true. \vskip 1mm

- The third step, known as the \emph{inductive step}, is to prove
that the given statement $P(i)$ (just assumed to be true) for any
one natural number $i$ implies that the given statement for the
next natural number $P(i+1)$ is true, i. e. to prove that
$P(i)\rightarrow P(i+1)$. \vskip 1mm

From these three steps, mathematical induction is the rule from
which we infer that the given statement $P(n)$ is established for
all natural numbers $n\geq k$.
\vskip 6mm

\section{The Lehmus-Steiner's theorem}

The Lehmus-Steiner's theorem states

\begin{thm}
If the straight lines bisecting the angles at the base of a triangle
and terminated by the opposite sides are equal, the triangle is
isosceles.
\end{thm}

This so called \emph{equal internal bisectors theorem} was
 communicated by Professor Lehmus (1780-1863) of Berlin
to Jacob Steiner (1796-1867) in the
year 1840 with a request for a pure geometrical
proof of it. The request was complied with at the time, but
Steiner's proof was not published till some years later. After
giving his proof, Steiner considers also the case when the angles
below the base are bisected; he generalizes the theorem somewhat;
found an external case where the theorem is not true; finally he
discusses the case of the spherical triangle. His solution by the
method of \emph{proof by Contraposition}  \cite{St} is considered
to be the most elementary one at that time.

Since then many mathematicians have published analytical and geometrical
solutions of this ``elementary'' theorem.

Does there exist a proof of this theorem which is \emph{direct}?
This problem was set in a Cambridge Examination Paper
in England around 1850. In 1853, the famous Bristish
mathematician James Joseph Sylvester (1814-1897)
intended to show that \emph{no direct proof can exist}, but he
was not very successful. Since then, there have been a
number of direct proofs published, but generally speaking they require some other results
which have not been proved directly.

A simple, constructive proof, based mainly on Euclid, Book III, is
given in \cite{W}.

McBride's paper \cite{M'B} contains a short history of the
theorem, a selection from the numerous other solutions that have
been published, some discussion of the logical points raised, and
a list of references to the extensive literature on the subject.

About the long history of this remarkable theorem see also \cite{M}.
\vskip 2mm

Here we propose two strictly direct proofs of Lehmus-Steiner's Theorem.

\vskip 4mm
\subsection{First proof of the Theorem of Lehmus-Steiner.}

\begin{figure}[h t b]
\epsfxsize=8cm \centerline{\epsfbox{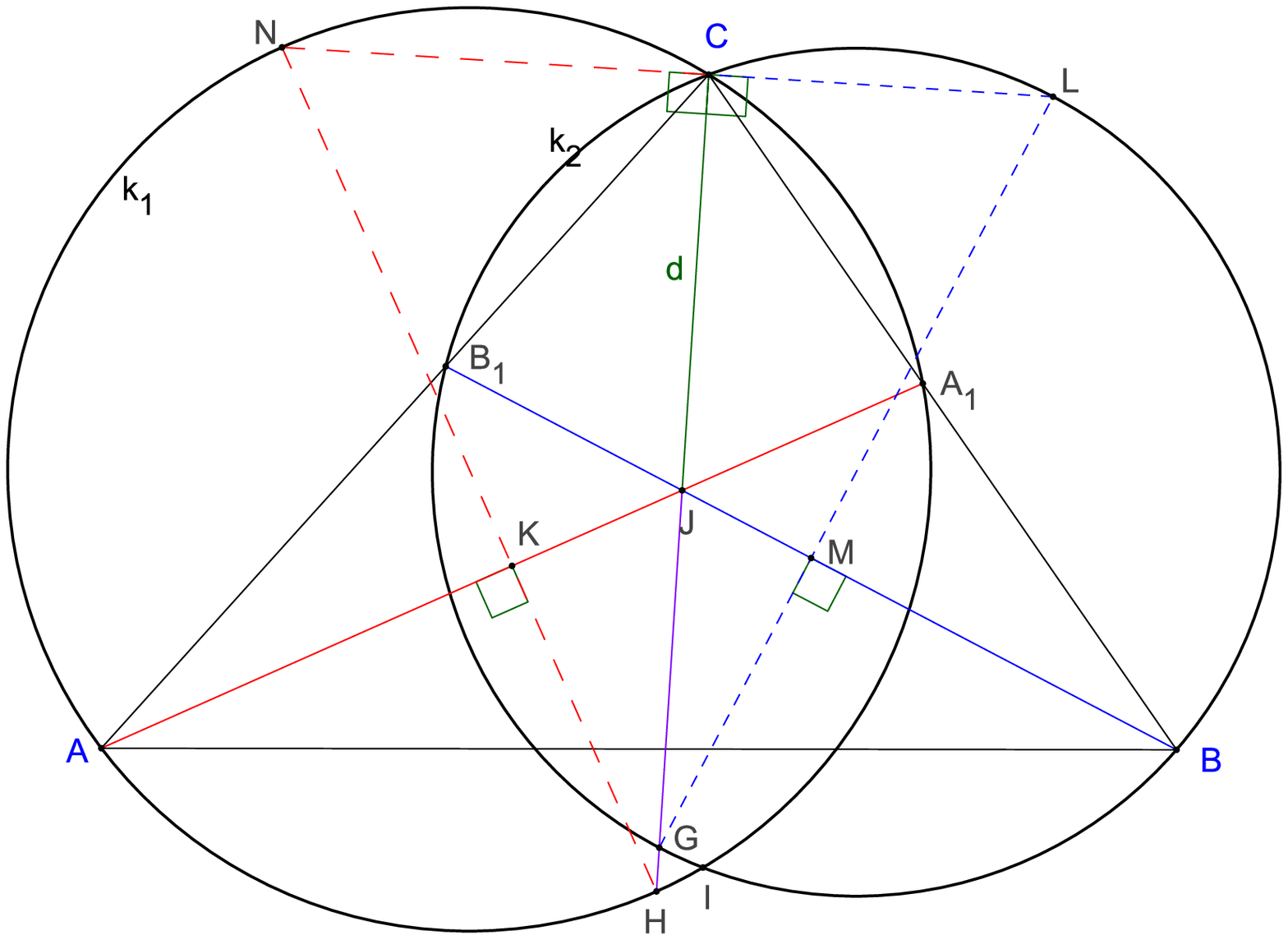}}
\end{figure}

Let  $AA_1\, (A_1\in BC)$  and $BB_1\, (B_1\in AC)$ be the
internal bisectors in $\,\triangle\, ABC$, $\,AA_1=BB_1$ and
$AA_1\cap BB_1=J$. Then $CJ$ is the internal bisector of $\angle
ACB$. We use the denotation $\;\gamma:=\angle ACJ =\angle BCJ.$

Let also $k_1$ be the circumscribing circle of $\triangle ACA_1$,
and $k_2$ the circumscribing circle of $\triangle BCB_1$ (fig. 1).

First we need the following
\begin{prop}The cut loci of points,  from which two equal segments appear under the same angle,
are equal arcs of congruent circles.
\end{prop}
\vskip 2mm

\emph{Proof of Proposition $3.2$.} Let us consider  $\triangle
ACA_1$ and $\triangle BC_1B_1$, where $\angle ACA_1=\angle
BC_1B_1= 2 \gamma$ and $AA_1=BB_1$. Let  $k_1$ with center $O_1$
be the circumscribing circle of $\triangle ACA_1$, and $k_2$ with
center $O_2$ the circumscribing circle of $\triangle BC_1B_1$
(fig. 2).

\begin{figure}[h t b]
\epsfxsize=8cm \centerline{\epsfbox{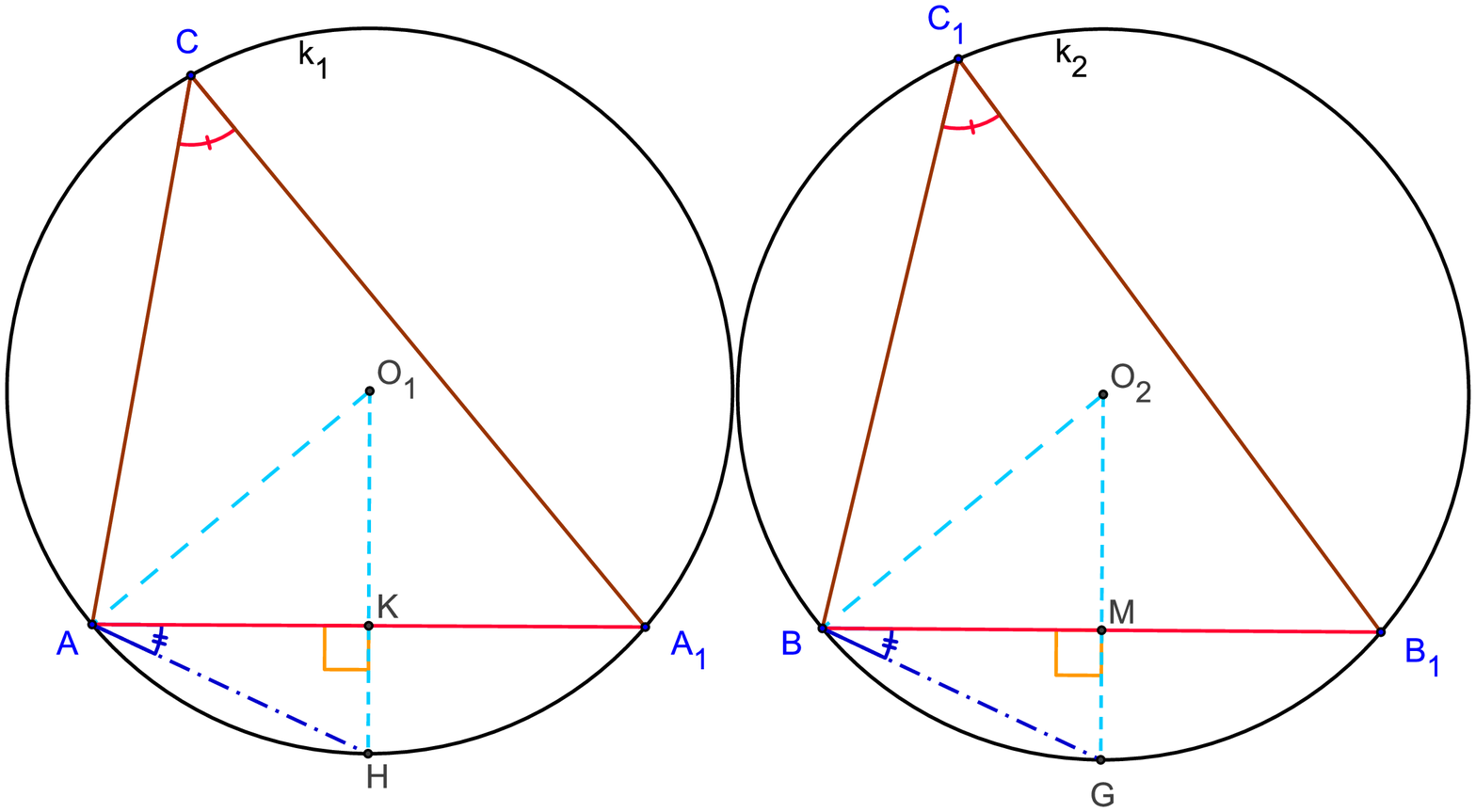}}
\end{figure}

The cut loci of points,  from which the equal segments $AA_1$ and
$BB_1$ appear under the same angle $2 \gamma$, are respectively
the arcs $\widehat{ACA_1}$ in $k_1$ and $\widehat{BC_1B_1}$ in $k_2$.

The perpendicular line $O_1K\; (K\in AA_1)$ from $O_1$ to the
chord $AA_1$ cuts the arc $\widehat{AA_1}$ in $k_1$ at its midpoint
$H$, the perpendicular line $O_2M\; (M\in BB_1)$ from $O_2$ to the
chord $BB_1$ cuts the arc $\widehat{BB_1}$ in $k_2$ at its midpoint
$G$.

The right angled triangles $\triangle AKH$ and $\triangle BMG$ are
congruent, because of $AK=BM$ (as a half of equal chords) and
$\angle KAH=\angle MBG=\gamma$. Hence, $AH=BG$ and $\angle AHK =
\angle BGM$.

Then, the isosceles triangles $\triangle AO_1H$ and $\triangle BO_2G$
are congruent and the circles $k_1$ and $k_2$ have equal radii.

This proves the assertion of the proposition.
$\hfill \square$
\vskip 4mm

Since the equal segments $AA_1$  and $BB_1$ in $\triangle ABC$
(fig. 1) appear under the same angle $2\gamma$ from $C$, the
circles $k_1$ and $k_2$ have equal radii (Proposition 3.2).
\vskip 1mm

Let now $CJ\cap k_1=H$ and $CJ\cap k_2=G$.

The points $H$ and $G$
lie on the same ray $CJ^{\overrightarrow{}}$. Since $CJ$ bisects the
angles $\angle ACA_1$ and $\angle BCB_1$, the point $H$ is
midpoint of the arc $\widehat{AA_1}$ in $k_1$, and the point $G$ is
midpoint of $\widehat{BB_1}$ in $k_2$.

Let $K$ be the midpoint of
the chord $AA_1$, $M$ be the midpoint of the chord $BB_1$,
$HK\cap k_1=N$ and $GM\cap k_2=L$.  Hence, the segments $HN$ and
$GL$ are diameters of the circles $k_1$ and $k_2$ respectively.
The triangles $\triangle CHN$ and $\triangle CGL$ are right angled
with right angles at the vertex $C$.

The quadrilateral $CJKN$ can be inscribed in a circle and it
follows that
$$|HK||HN|=|HJ||HC|. \leqno (1)$$

The quadrilateral $CJML$ can be inscribed in a circle and it
follows that
$$|GM||GL|=|GJ||GC|. \leqno (2)$$
\vskip 2mm

\begin{rem}
The equalities (1) and (2) are also a consequence of the
similarities

$\;\triangle HKJ\sim \triangle HCN\;$  and $\;\triangle GMJ\sim
\triangle GCL$.
\end{rem}
\vskip 2mm

Since the circles $k_1$ and $k_2$ have equal radii and the chords
$AA_1$  and $BB_1$ are equal, then $HK=GM$ and $HN=GL$. If we put
$\;d=|CJ|>0$, $\;x=|HJ|>0$, $\;y=|GJ|>0$, then $|HC|=x+d$ and
$|GC|=y+d$.

The left hand sides of equalities (1) and (2) are equal, so
are their right hand sides. Hence
$$x(x+d)=y(y+d)\;\Leftrightarrow\; (x-y)(x+y+d)=0. \leqno (3)$$

Since  $\,x+y+d\neq 0$, equality (3) is equivalent to the equality
$$ x-y=0\,.\,\frac{1}{x+y+d}\,=0,$$

which directly implies  $\;x=y$.
\vskip 2mm

\begin{rem} If we denote the equal positive left hand sides of equalities (1) and (2) by $a^2$,
we get respectively the quadratic equations
$$\begin{array}{ll}
x^2+dx-a^2=0\; &\Leftrightarrow\;\displaystyle{\left(x+\frac{d}{2}\right)^2-
\left(\frac{\sqrt{4\,a^2+d^2}}{2}\right)^2=0}\\
[6mm]
\;&\Leftrightarrow\;\displaystyle{\left(x+\frac{\sqrt{4\,a^2+d^2}+d}{2}\right)
\left(x-\frac{\sqrt{4\,a^2+d^2}-d}{2}\right)=0}\\
[6mm]
\;&\Leftrightarrow\;\displaystyle{x-\frac{\sqrt{4\,a^2+d^2}-d}{2}=
0\,.\,\left(x+\frac{\sqrt{4\,a^2+d^2}+d}{2}\right)^{-1}=0,}
\end{array}$$

and analogously

$$\quad y^2+dy-a^2=0\;\Leftrightarrow\;
y-\frac{\sqrt{4\,a^2+d^2}-d}{2}=
0\,.\,\left(y+\frac{\sqrt{4\,a^2+d^2}+d}{2}\right)^{-1}=0$$

with the same solution

$$x=y=\frac{1}{2}\left(\sqrt{4a^2+d^2}-d\right).$$
\end{rem}
\vskip 4mm

Hence, the points $H$ and $G$, which lie on the same ray, coincide and $CG$ is the common chord of
the circles $k_1$ and $k_2$.
\vskip 2mm

As a consequence of the fulfilled conditions
\begin{itemize}
\item[-] $CG$ is a common side,
\item[-] $\angle ACG=\angle BCG$
($CG$ is the bisector of $\angle ACB$),
\item[-] $\angle CAG=\angle CBG$
($CG$ is the common chord of two circles with equal radii, hence
$\widehat{CA_1G}=\widehat{CB_1G}$),
\end{itemize}
the triangles $\triangle AGC$ and $\triangle BGC$ are congruent
(fig. 3)

\begin{figure}[h t b]
\epsfxsize=8cm \centerline{\epsfbox{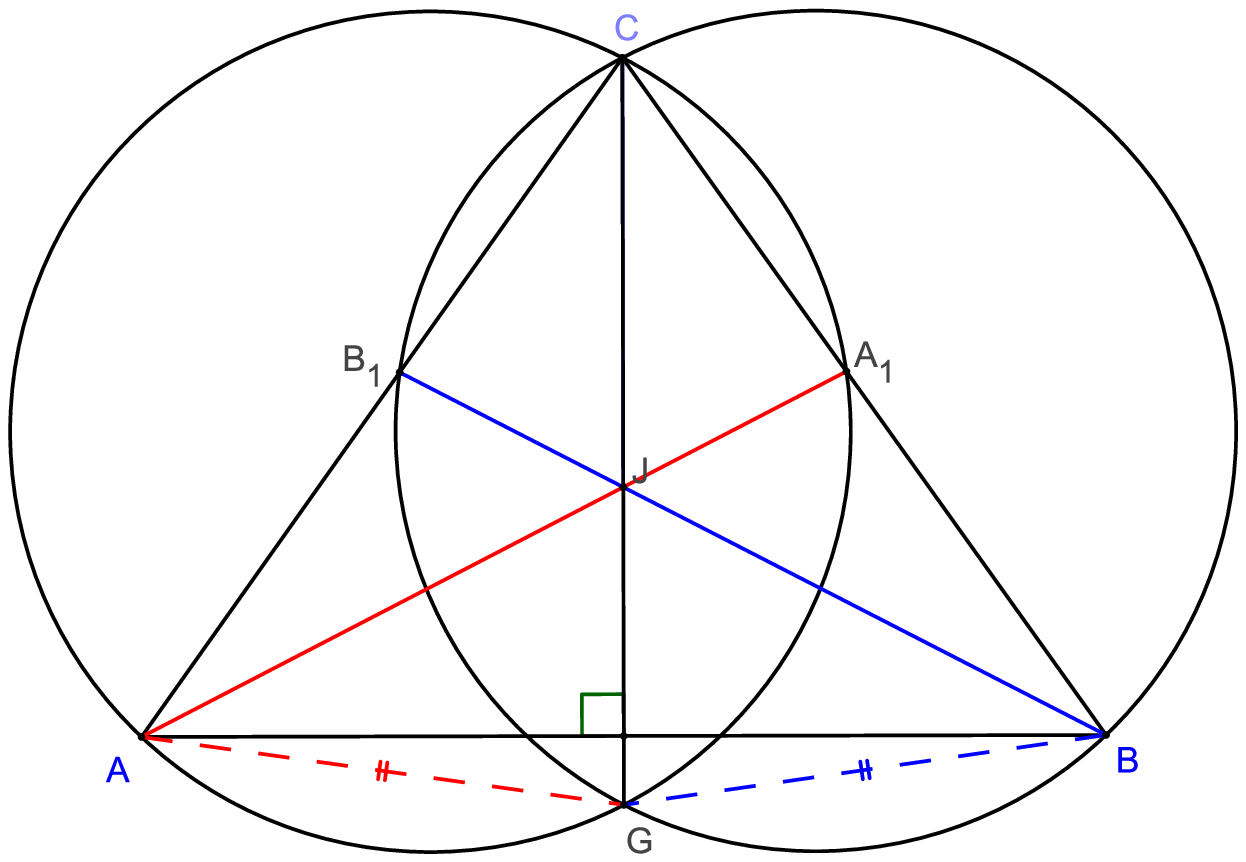}}
\end{figure}
\vskip 1mm

Thus, $CA=CB$ and  $\triangle ABC$ is isosceles.

The direct proof of Theorem 3.1 is complete.
$\hfill\square$ \vskip 4mm

\begin{rem}
In this proof, the condition the segments  $AA_1$ and $BB_1$
are internal bisectors of the angles based at $AB$ in $\triangle ABC$ is
not necessary.

It is only of importance that they are equal by length
cevians and their intersection point lies on the bisector of $\angle ACB$.

We recall that a \emph{cevian} is a line segment which joins a vertex of a triangle
with a point on the opposite side (or its extension).
\vskip 2mm

In fact we  proved directly the following
\begin{thm}
If in a $\triangle ABC$ the segments $AA_1\, (A_1\in BC)$ and
$BB_1\, (B_1\in AC)$ cut at a point on the bisector of  $\angle
ACB$ and are equal by length then   $\triangle \,ABC$ is
isosceles.
\end{thm}
\end{rem}
\vskip 4mm

\subsection{Second proof of the Theorem of Lehmus-Steiner.} The idea for this proof comes
from Problem 2.1-16 in \cite{E}: \emph{Find a direct proof of
Lehmus-Steiner's theorem as a consequence of Stewart's theorem.}
\vskip 2mm

We need the notion \emph{algebraic measure} (\emph{relative measure}) of a line segment.

On any straight line there are two (opposite to each other)
directions. The \emph{axis} is a couple of a straight line and a fixed (positive)
direction on it.

Let $g^+$ denotes any axis. For any non zero line segment
$MN$ on $g^+$ we can define its \emph{relative} (\emph{algebraic}) \emph{measure} by
$\overline{MN}=\varepsilon |MN|$, where $\varepsilon =+1$ in case
$\overrightarrow{MN}$ has the same direction as $g^+$, and
$\varepsilon =-1$ in case $\overrightarrow{MN}$ has the opposite
direction with respect to  $g^+$.
\vskip 2mm

 \emph{Stewart's theorem} yields a relation between the lengths of the sides of a
 triangle and the length of a cevian.
\vskip 2mm

 Let in $\triangle \,ABC$ the line segment $CP,\, P\in AB,$ be a cevian (more general, let
 $\{C; A, B, P\}$ be a quadruple of points such that $A, B, P$ are collinear).
\begin{thm} (Theorem of Stewart)
If $A, B$,  $P$ are three collinear points and $C$ is any
point then
$$|CA|^2\cdot \overline{BP}+|CB|^2\cdot \overline{PA}+|CP|^2\cdot \overline{AB}
+ \overline{BP}\cdot \overline{PA}\cdot \overline{AB} =0.$$
\end{thm}

\begin{rem}
Using the Theorem of Pythagoras, the proof of Steward's theorem  is a simply verification.
\end{rem}
\vskip 3mm

In what follows we prove the \emph{equal internal bisectors theorem} in the following formulation.
\begin{thm}
The straight lines bisecting the angles at the base of a triangle
and terminated by the opposite sides are equal if and only if the triangle is
isosceles.
\end{thm}

Let   $AA_1\, (A_1\in BC)$  and $BB_1\, (B_1\in AC)$ be
respectively the internal bisectors of $\angle CAB$ and $\angle
CBA$ in a triangle $ ABC$ (fig. 4).

\begin{figure}[h t b]
\epsfxsize=8cm \centerline{\epsfbox{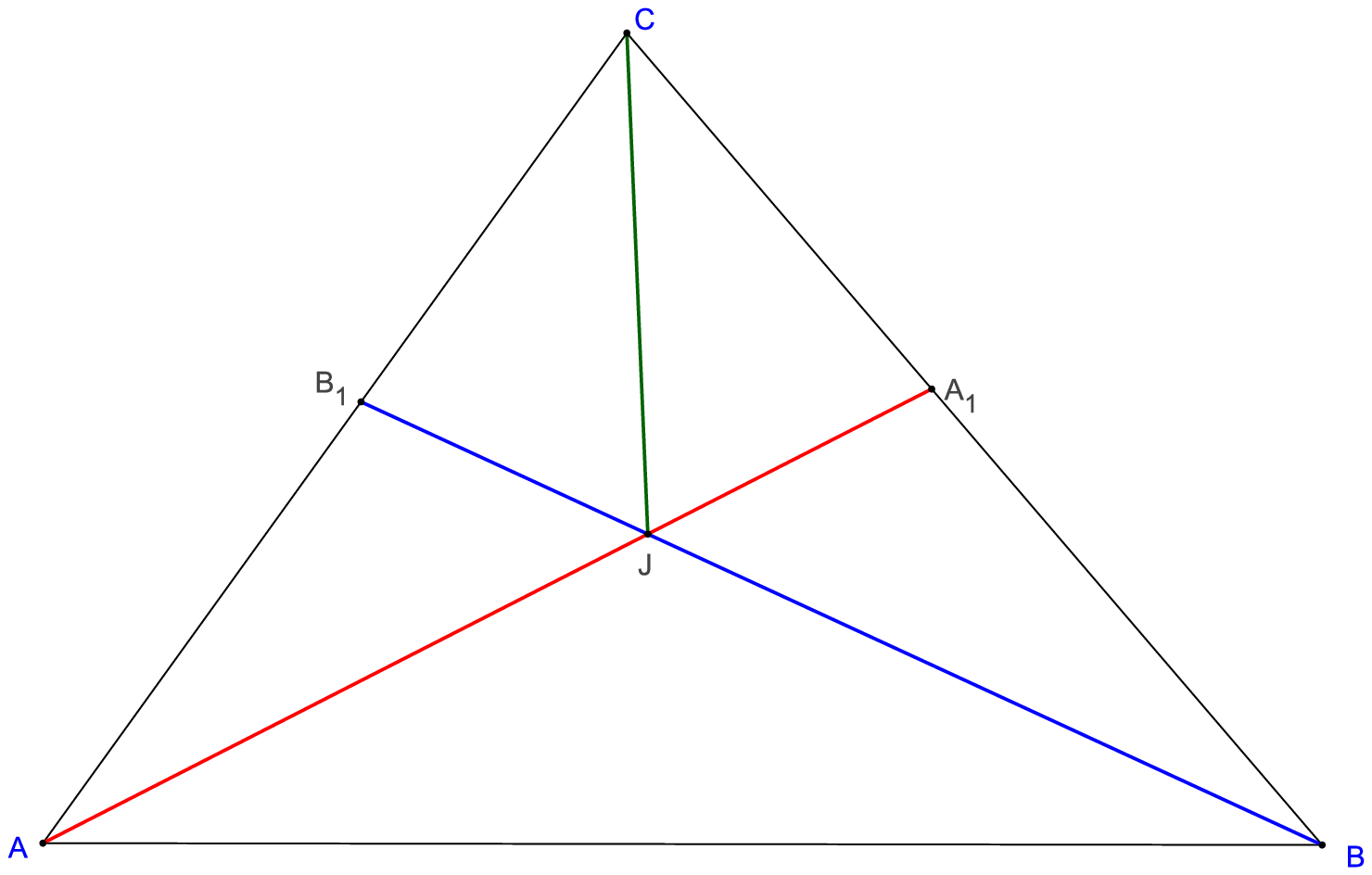}}
\end{figure}

Since the triples $\{B, A_1, C\}$ and $\{A, B_1, C,\}$ consist of
collinear points there exist integers $\alpha$ and $\beta$ such
that
$$\begin{array}{ll}
\overline{BA_1}=\alpha\,\overline{BC},\quad & \overline{A_1C}=(1-\alpha)\,\overline{BC},
\quad 0<\alpha< 1;\\
[2mm]
\overline{AB_1}=\beta\,\overline{AC},\quad & \overline{B_1C}=(1-\beta)\,\overline{AC},
\quad 0<\beta< 1.
\end{array} \leqno(4)$$

Using the fact that $AA_1\, (A_1\in BC)$  and $BB_1\, (B_1\in AC)$
are the internal bisectors of $\angle CAB$ and $\angle CBA$ in a
triangle $ ABC$, i. e. that
$$\frac{\overline{CA_1}}{\overline{A_1B}}=\frac{|CA|}{|AB|},\quad
\frac{\overline{CB_1}}{\overline{B_1A}}=\frac{|CB|}{|BA|},$$
 from  relations (4) we compute

$$\begin{array}{ll}
\displaystyle{\alpha=\frac{|AB|}{|AB|+|AC|}}\,,\quad & \displaystyle{ 1-\alpha=\frac{|AC|}{|AB|+|AC|}\,,}\\
[6mm]
\displaystyle{\beta=\frac{|AB|}{|AB|+|BC|}}\,,\quad & \displaystyle{1-\beta=\frac{|BC|}{|AB|+|BC|}\,.}
\end{array} \leqno(5)$$
\vskip 2mm

Applying Stewart's theorem for the quadruple $\{A; B, A_1, C\}$
$$|AB|^2.\overline{A_1C}+|AA_1|^2.\overline{CB}+|AC|^2.\overline{BA_1}+
\overline{A_1C}.\overline{CB}.\overline{BA_1}=0,$$
and for the quadruple $\{B; A, B_1, C\}$
$$|BA|^2.\overline{B_1C}+|BB_1|^2.\overline{CA}+|BC|^2.\overline{AB_1}+
\overline{B_1C}.\overline{CA}.\overline{AB_1}=0,$$
from (4) and (5) we get
$$|AA_1|^2=\frac{|AB||AC|}{(|AB|+|AC|)^2}\,\{(|AB|+|AC|)^2-|BC|^2\},$$
\vskip 4mm

$$|BB_1|^2=\frac{|AB||BC|}{(|AB|+|BC|)^2}\,\{(|AB|+|BC|)^2-|AC|^2\},$$
and finally

$$\begin{array}{l}
\displaystyle{(|AA_1|-|BB_1|)\,\frac{(|AA_1|+|BB_1|)}{|AB|}=}\\
[8mm]
\displaystyle{(|AC|-|BC|)\left\{1+\frac{|AC||BC|(|AB|^2+|AC|^2+|BC|^2
+2|AB|(|AC|+|BC|)+|AC||BC|)}{(|AB|+|BC|)^2(|AB|+|AC|)^2}\right\}.}
\end{array}$$
\vskip 3mm

Using the denotations
$$X:= \frac{(|AA_1|+|BB_1|)}{|AB|}$$

\noindent
and
$$Y:=\left\{1+\frac{|AC||BC|(|AB|^2+|AC|^2+|BC|^2
+2|AB|(|AC|+|BC|)+|AC||BC|)}{(|AB|+|BC|)^2(|AB|+|AC|)^2}\right\},$$

\noindent
we rewrite the last equation in the form

$$(|AA_1|-|BB_1|)\,X=(|AC|-|BC|)\,Y.\leqno (6)$$

Since $X\neq 0$ and $Y\neq 0$, equation (6) is equivalent to the equations
$$\left(|AA_1|-|BB_1|\right)\,\frac{X}{Y}=|AC|-|BC|\leqno (7)$$
and
$$\left(|AC|-|BC|\right)\,\frac{Y}{X}=|AA_1|-|BB_1|.\leqno (8)$$
\vskip 4mm

Thus, if $|\;AA_1|=|BB_1|$, it results from equality (7) directly that $\;|AC|=|BC|$, i. e.
$$|AA_1|=|BB_1|\;\Rightarrow\;|AC|=|BC|.$$
\vskip 4mm

If $\;|AC|=|BC|$, it results from equality (8) directly that $|AA_1|=|BB_1|$, i. e.
$$|AC|=|BC|\;\Rightarrow\;|AA_1|=|BB_1|.$$
\vskip 4mm

Hence,
$$|AA_1|=|BB_1|\quad \Leftrightarrow\quad |AC|=|BC|,$$
which completes this direct proof of Lehmus-Steiner's theorem.
\hfill{$\square$}
\vskip 2mm

\begin{rem}
\begin{itemize}
\item In this proof, the condition the segments  $AA_1$ and $BB_1$
are internal bisectors of the angles based at $AB$ in $\triangle ABC$ is
necessary.
\vskip 2mm

\item Using equalities (5) we compute
$$\alpha - \beta= \frac{|AB|}{(|AB|+|AC|)(|AB|+|BC|)}\,(|BC|-|AC|)$$
and obtain
$$A_1B_1\,\parallel \,
AB\quad\Leftrightarrow\quad\alpha=\beta\quad\Leftrightarrow\quad |AC|=|BC|.$$
\end{itemize}
\vskip 4mm

The following statement is easily to be proved directly.

\begin{prop}
Let   $AA_1\, (A_1\in BC)$  and $BB_1\, (B_1\in AC)$ be
respectively the internal bisectors of $\angle CAB$ and $\angle
CBA$ in  $ \triangle \,ABC$. Prove that
$\triangle\, ABC$ is isosceles if and only if $A_1B_1\,\parallel\,AB$.
\end{prop}

\emph{Proof}. Let $AA_1\, (A_1\in BC)$  and $BB_1\, (B_1\in AC)$ be
respectively the internal bisectors of $\angle CAB$ and $\angle
CBA$ in  $ \triangle \,ABC$.

\begin{itemize}
\item[-] Let $A_1B_1\,\parallel\,AB$ (fig. 5).

It follows that $\, \triangle\, AA_1B_1\,$ and $\, \triangle\, BB_1A_1\,$
are isosceles and the quadrilateral $\,ABA_1B_1\,$ is a trapezium with  $\,|AB_1|=|BA_1| \,(=|A_1B_1|)$.

Hence, $\triangle\, ABC$ is isosceles.
\vskip 2mm

\begin{figure}[h t b]
\epsfxsize=8cm \centerline{\epsfbox{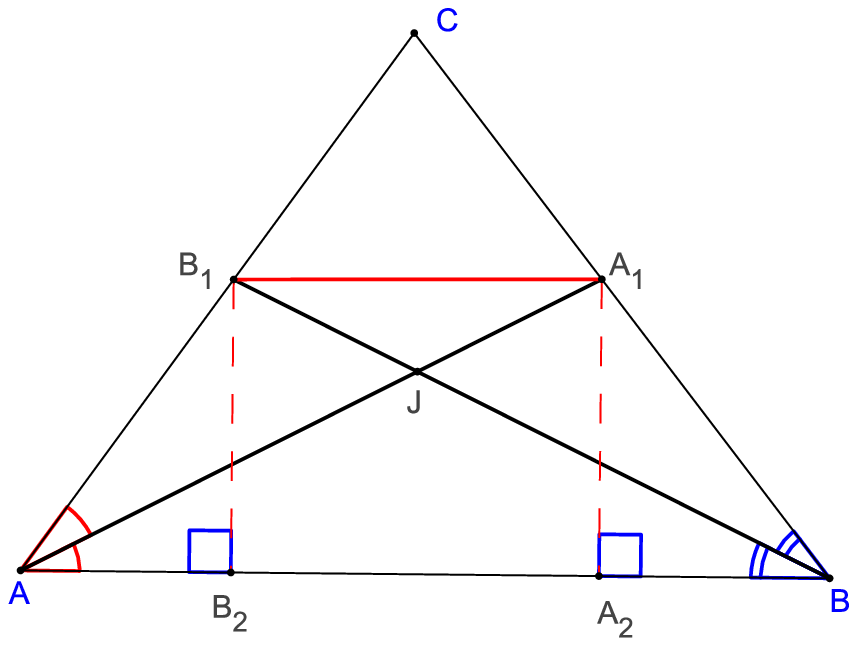}}
\end{figure}

\item[-] Let now $\triangle\, ABC$ be isosceles and $\,B_1B_2\perp AB$ $\,(B_2\in AB)$,
$A_1A_2\perp AB$ $\,(A_2\in AB)$.

Since $\,\triangle\,AA_1B\cong\triangle\, BB_1A$ (fig. 5), then $\,|AA_1|=|BB_1|$.

Hence, $\,\triangle\,AA_1A_2\cong\triangle\, BB_1B_2$, $|A_1A_2|=|B_1B_2|\,$ and $\, A_1B_1\parallel AB$.
\end{itemize}
\hfill{$\square$}
\vskip 2mm

In view of this Proposition we can reformulate the Theorem of Lehmus-Steiner in the form:

\begin{thm}
Let   $AA_1\, (A_1\in BC)$  and $BB_1\, (B_1\in AC)$ be
respectively the internal bisectors of $\angle CAB$ and $\angle
CBA$ in  $ \triangle \,ABC$. Prove that if $|AA_1|=|BB_1|$, then
$A_1B_1\,\parallel\,AB$.
\end{thm}
\end{rem}

\vskip 4mm
Acknowledgements. The first author is partially supported by Sofia University Grant 59/2014.
The second author is partially supported by Sofia University Grant 126/2014.

\end{document}